\documentclass[a4paper,10pt]{amsart}

\usepackage[english]{babel}

\usepackage[utf8]{inputenc}

\usepackage{amsmath,amsfonts,amssymb,latexsym,graphics,amsthm}
\usepackage[bookmarks=false]{hyperref}
\usepackage[dvips]{color}
\usepackage[dvips]{graphicx}
\usepackage{pstricks}
\usepackage{url}

\renewenvironment{itemize}{\begin{list}{\labelitemi}{\leftmargin=1.5em}}{\end{list}}
\renewcommand{\labelitemi}{$\bullet$}

\newcommand{\ket}[1]{\ensuremath{|#1\rangle}}
\newcommand{\bra}[1]{\ensuremath{\langle #1|}}
\newcommand{\braket}[2]{\ensuremath{\langle #1|#2 \rangle}}
\newcommand{\Le}{\hbox{\rotatedown{$\Gamma$}}}

\title{A $q$-enumeration of alternating permutations}
\date{\today}
\author{Matthieu Josuat-Vergès}
  \thanks{Partially supported by the ANR Jeune Chercheur IComb. A short version of this work
  was presented at the conference Permutation Patterns 2009.}
\address{LRI, CNRS and Université Paris-Sud, 
  Bâtiment 490, 91405 Orsay, FRANCE}
\email{josuat@lri.fr}

\newtheorem{thm}{Theorem}[section]

\newtheorem{lem}[thm]{Lemma}
\newtheorem{prop}[thm]{Proposition}
\newtheorem{defn}{Definition}

\theoremstyle{definition}
\newtheorem{rem}[thm]{Remark}

\begin{document}

\begin{abstract}
A classical result of Euler states that the tangent numbers are an alternating
sum of Eulerian numbers. A dual result of Roselle states that the secant numbers
can be obtained by a signed enumeration of derangements. We show that both 
identities can be refined with the following statistics: the number of crossings
in permutations and derangements, and the number of patterns 31-2 in alternating
permutations.

Using previous results of Corteel, Rubey, Prellberg, and the author, we derive closed 
formulas for both $q$-tangent and $q$-secant numbers. There are two different methods
to obtain these formulas: one with permutation tableaux and one with weighted Motzkin 
paths (Laguerre histories).
\end{abstract}

\maketitle

\section{Introduction}

The classical Euler numbers $E_n$ are given by the Taylor expansion of the tangent and 
secant functions: 
\[   tan(x) + sec(x) = \sum_{n=0}^{\infty} E_n\frac{x^n}{n!}.  \]
Since the tangent is an odd function and the secant is an even function, the integers 
$E_{2n+1}$ are called the {\it tangent numbers} and the integers $E_{2n}$ are called the 
{\it secant numbers}. More than a century ago, Désiré André \cite{DA} showed that 
$E_n$ is the number of {\it alternating permutations} in $\mathfrak{S}_n$.

Several $q$-analogs of these numbers have been studied, mainly by Andrews, Foata, Gessel,
Han \cite{AF,AG,FH,HZR}. The ones we consider here appeared in an article
of Han, Randrianarivony, Zeng \cite{HZR}.

\begin{defn} The $q$-tangent numbers $E_{2n+1}(q)$ are defined by:
\begin{equation}
\sum_{n=0}^{\infty} E_{2n+1}(q) x^n = 
\cfrac{1}{1-\cfrac{[1]_q[2]_qx}{1-\cfrac{[2]_q[3]_qx}{1-\cfrac{[3]_q[4]_qx}{\ddots}}}},
\end{equation}
where $[n]_q=\frac{1-q^n}{1-q}$. And the $q$-secant numbers $E_{2n}(q)$ are defined by:
\begin{equation}
\sum_{n=0}^{\infty} E_{2n}(q) x^n = 
\cfrac{1}{1-\cfrac{[1]_q^2x}{1-\cfrac{[2]_q^2x}{1-\cfrac{[3]_q^2x}{\ddots}}}}.
\end{equation}
\end{defn}

These continued fractions are the natural $q$-analogs of the ones appearing in
Flajolet's celebrated article \cite{Fla}. The methods in this reference gives 
a combinatorial interpretation of $E_n(q)$ in terms of weighted Dyck paths. 
Indeed, if $\delta\in\{0,1\}$ then $E_{2n+\delta}(q)$ is the sum of weights of 
Dyck paths of length $2n$ such that:
\begin{itemize}
\item the weight of a step $\nearrow$ starting at height $h$ is $q^i$ for some 
  $i\in\left\{0,\dots,h\right\}$,
\item the weight of a step $\searrow$ starting at height $h$ is $q^i$ for some 
  $i\in\left\{0,\dots,h-1+\delta\right\}$,
\end{itemize}
and the weight of the path is the product of the weights of its steps.
The first values are $E_0(q)=E_1(q)=E_2(q)=1$, $E_3(q)=1+q$, $E_4(q)=2+2q+q^2$, 
$E_5(q)=2+5q+5q^2+3q^3+q^4$.

\medskip

Another combinatorial interpretation is that
\begin{equation} \label{combint}
E_n(q) = \sum_{\sigma\in \mathfrak{A}_n} q^{31\hbox{-}2(\sigma)},
\end{equation}
where $31$-$2(\sigma)$ is the number of occurrences of the generalized pattern $31$-$2$,
and $\mathfrak{A}_n$ is the set of alternating permutations of size $n$.
This was proved bijectively by D. Chebikin \cite[Theorem 9.5]{Che}. Another bijective proof 
can be obtained with the refined Françon-Viennot bijection given by Corteel \cite{Co}. 
For more details about this see Remarks \ref{rem-even} and \ref{rem-odd} below.

\medskip

Let $wex(\sigma)$ be the number of weak exceedances of $\sigma$.
In \cite{FS} Foata and Sch\"utzenberger show that the specialization of the Eulerian 
polynomial $\sum_{\sigma\in\mathfrak{S}_n} x^{wex(\sigma)}$ at $x=-1$ is 0 when $n$ is 
even and $(-1)^{\frac{n-1}2}E_{n}$ when $n$ is odd. This was first obtained by Euler
\cite{Euler}.
Let $\mathfrak{D}_n\subset\mathfrak{S}_n$ be the set of derangements, then another result 
is that $\sum_{\sigma\in\mathfrak{D}_n} x^{wex(\sigma)}$ specialized at $x=-1$ is 
$(-1)^{\frac n2}E_{n}$ when $n$ is even and 0 when $n$ is odd. This was first obtained by
Roselle \cite{DPR}.

\medskip

Refinements of these two results have been given by Foata and Han \cite{FH}. The involved
statistic in permutations and derangements is the the major index, and the involved statistic 
in alternating permutations is the number of inversions. In this article we give different 
refinements of the same equalities. The relevant statistic is the number of {\it crossings} 
of a permutation, denoted by  $cr(\sigma)$ and introduced in \cite{Co}. A crossing of 
$\sigma\in\mathfrak{S}_n$ is a couple $(i,j)$ such that either $i<j\leq\sigma(i)<\sigma(j)$ 
or $\sigma(i)<\sigma(j)<i<j$.

\begin{thm} \label{th1}
\begin{equation} \label{eq1}
\sum_{\sigma\in\mathfrak{S}_n} (-1)^{wex(\sigma)} q^{cr(\sigma)} = 
\begin{cases} 0 & \mbox{if } n\mbox{ is even,} \\ (-1)^{\frac{n-1}2}E_{n}(q) & \mbox{if } 
n\mbox{ is odd.} \end{cases}
\end{equation}
\end{thm}

\begin{thm} \label{th2}
\begin{equation} \label{eq2}
\sum_{\sigma\in\mathfrak{D}_n} \left(-\tfrac 1q\right)^{wex(\sigma)} q^{cr(\sigma)} = 
\begin{cases} \left(-\tfrac 1q\right)^{\frac n2} E_{n}(q) & \mbox{if } n\mbox{ is even,}
\\ 0 & \mbox{if }
n\mbox{ is odd.} \end{cases}
\end{equation}
\end{thm}

The bistatistic $(wex,cr)$ is equidistributed with $(asc,31\hbox{-}2)$ in $\mathfrak{S}_n$ 
(but not in $\mathfrak{D}_n$). It is also equidistributed with the number of rows and
number of superfluous 1s in {\it permutation tableaux} \cite{SW}. From the latter description 
we see that the left-hand side of (\ref{eq1}) is an alternating sum of the $q$-Eulerian numbers
$\hat E_{k,n}(q)$ defined in \cite{LW}. In this reference, L. Williams gives multivariate 
generating functions for permutation tableaux, and a formula for $q$-Eulerian numbers that reduces 
to the classical formula of eulerian numbers when $q=1$:
\[ \hat E_{k,n}(q) = \sum_{i=0}^{k-1}(-1)^i[k-i]^nq^{ki-k^2} \left( \tbinom n i q^{k-i} 
    + \tbinom n{i-1} \right). 
\]
Another expression was conjectured by Corteel and Rubey
and proved in \cite{CJPR,MJV}:
\begin{equation}  \label{An}   \begin{split}
   \sum_{k=0}^n y^k \hat E_{k,n}(q) = 
   \tfrac 1{(1-q)^n}
  \sum\limits_{k=0}^n (-1)^k
    \Bigg( \sum\limits_{j=0}^{n-k} y^j\Big( \tbinom{n}{j}\tbinom{n}{j+k} -
    \tbinom{n}{j-1}\tbinom{n}{j+k+1}\Big) \Bigg) \\ \times
  \left( \sum\limits_{i=0}^k y^iq^{i(k+1-i)} \right).
\end{split}\end{equation}
There are two different proofs, that use respectively the combinatorics of Laguerre histories
and permutation tableaux. With similar methods, we can obtain closed formulas for $E_n(q)$.

\begin{thm} \label{th3}
\begin{equation} \label{En}
E_{2n+1}(q) = \frac 1{(1-q)^{2n+1}} \sum_{k=0}^n \left( \tbinom{2n+1}{n-k}-
\tbinom{2n+1}{n-k-1} \right)
\sum_{i=0}^{2k+1} (-1)^{i+k} q^{i(2k+2-i)}.
\end{equation}
\end{thm}

\begin{thm} \label{th4}
\begin{equation}
E_{2n}(q) = \frac 1{(1-q)^{2n}} \sum_{k=0}^n \left( \tbinom{2n}{n-k}-
\tbinom{2n}{n-k-1} \right)
\sum_{i=0}^{2k} (-1)^{i+k} q^{i(2k-i)+k}.
\end{equation}
\end{thm}

We can notice similarities with the Touchard-Riordan formula \cite{JGP}, which gives
the distribution of crossings in the set $\mathfrak{I}_{2n}\subset\mathfrak{S}_{2n}$
of fixed-point free involutions. If $\sigma\in\mathfrak{I}_{2n}$, then 
$\frac 12 cr(\sigma)$ is the number of couples $(i,j)$ such that $i<j<\sigma(i)<\sigma(j)$,
so that we recover a more classical definition of crossings. The formula is:
\begin{equation} \label{tou_rio}
\sum_{I \in \mathfrak{I}(2n)} q^{\frac 12cr(I)} = \frac 1 {(1-q)^n}
\sum_{k=0}^n \left( \tbinom{2n}{n-k} - \tbinom{2n}{n-k-1} \right) (-1)^k
q^{\frac {k(k+1)}2}.
\end{equation}
Note that the left-hand side of (\ref{tou_rio}) is also the sum of weights of Dyck 
paths such that the weight of a step $\nearrow$ is always 1, and the weight of a 
step $\searrow$ starting at height $h$ is $q^i$ for some $i\in\left\{0,\dots,h-1\right\}$. 
So this is a subset of the paths for $E_{2n}(q)$.

\bigskip

This article is organized as follows. Section \ref{sec2} contains preliminaries
that are mostly contained in the references. Section \ref{sec3} contains the proofs
of Theorems \ref{th1} and  \ref{th3}. Section \ref{sec4} contains the proofs
of Theorems \ref{th2} and  \ref{th4}. Section \ref{sec5} contains alternative 
proofs of Theorems \ref{th3} and  \ref{th4}.

\section*{Acknowledgement}

I thank my advisor Sylvie Corteel for her help, advice, support and kindness.
I thank Dominique Foata, Jiang Zeng, Guo-Niu Han, and Xavier Viennot for their insightful
suggestions.

\section{Definitions, conventions, preliminaries}
\label{sec2}

\subsection{Permutations and Laguerre histories}
Through this article, we use the convention that $\sigma(0)=0$ and $\sigma(n+1)=n+1$
if $\sigma\in\mathfrak{S}_n$. Let $i\in\{1,\dots,n\}$, it is a {\it weak exceedance} of $\sigma$ 
if $\sigma(i)\geq i$, an {\it ascent} of $\sigma$ if $\sigma(i)<\sigma(i+1)$ (note that $n$ is 
always an ascent with our conventions). We denote by $wex(\sigma)$ and $asc(\sigma)$ the total 
number of weak exceedances and ascents in $\sigma$.
A {\it descent} of $\sigma$ is a $i$ such that $\sigma(i)>\sigma(i+1)$.

A {\it Laguerre history} (or ``histoire de Laguerre'') of size $n$ is a weighted Motzkin path of $n$ 
steps such that:
\begin{itemize}
\item the weight of a step $\nearrow$ starting at height $h$ is $yq^i$ for some
  $i\in\left\{0,\dots,h\right\}$,
\item the weight of a step $\rightarrow$ starting at height $h$ is either $yq^i$ for
  some $i\in\left\{0,\dots,h\right\}$ or $q^i$ for some $i\in\left\{0,\dots,h-1\right\}$,
\item the weight of a step $\searrow$ starting at height $h$ is $q^i$ for some
  $i\in\left\{0,\dots,h-1\right\}$.
\end{itemize}

They are in bijection with other kinds of Motzkin paths with one fewer step.
A {\it large Laguerre history} of size $n$ is a weighted Motzkin path of $n-1$ steps
such that:
\begin{itemize}
\item the weight of a step $\nearrow$ starting at height $h$ is $yq^i$ for some
  $i\in\left\{0,\dots,h\right\}$,
\item the weight of a step $\rightarrow$ starting at height $h$ is either $yq^i$  
  or $q^i$ for some $i\in\left\{0,\dots,h\right\}$,
\item the weight of a step $\searrow$ starting at height $h$ is $q^i$ for some
  $i\in\left\{0,\dots,h\right\}$.
\end{itemize}

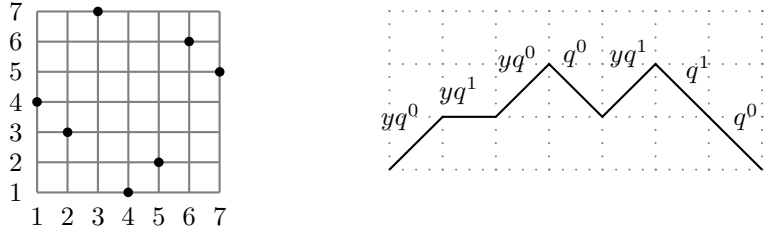
\begin{figure}[h!tp]\psset{unit=4mm}
\begin{pspicture}(0,0)(7,7)
\rput(0.3,1){1}
\rput(0.3,2){2}
\rput(0.3,3){3}
\rput(0.3,4){4}
\rput(0.3,5){5}
\rput(0.3,6){6}
\rput(0.3,7){7}
\rput(1,0.2){1}
\rput(2,0.2){2}
\rput(3,0.2){3}
\rput(4,0.2){4}
\rput(5,0.2){5}
\rput(6,0.2){6}
\rput(7,0.2){7}
\psline[linecolor=gray](1,1)(1,7)\psline[linecolor=gray](1,1)(7,1)
\psline[linecolor=gray](2,1)(2,7)\psline[linecolor=gray](1,2)(7,2)
\psline[linecolor=gray](3,1)(3,7)\psline[linecolor=gray](1,3)(7,3)
\psline[linecolor=gray](4,1)(4,7)\psline[linecolor=gray](1,4)(7,4)
\psline[linecolor=gray](5,1)(5,7)\psline[linecolor=gray](1,5)(7,5)
\psline[linecolor=gray](6,1)(6,7)\psline[linecolor=gray](1,6)(7,6)
\psline[linecolor=gray](7,1)(7,7)\psline[linecolor=gray](1,7)(7,7)
\psdots(1,4)(2,3)(3,7)(4,1)(5,2)(6,6)(7,5)
\end{pspicture}  \hspace{2cm} \psset{unit=7mm}
\begin{pspicture}(0,-1)(7,3)
\psgrid[gridcolor=gray,griddots=4,subgriddiv=0,gridlabels=0](0,0)(7,3)
\psline(0,0)(1,1)(2,1)(3,2)(4,1)(5,2)(6,1)(7,0)
\rput(0.2,1.0){$yq^0$}
\rput(1.3,1.5){$yq^1$}
\rput(2.4,2.1){$yq^0$}
\rput(3.5,2.2){$q^0$}
\rput(4.5,2.2){$yq^1$}
\rput(5.8,1.9){$q^1$}
\rput(6.7,1){$q^0$}
\end{pspicture}
\caption{Example of the permutation $4371265$ and its image by the Françon-Viennot 
bijection. \label{fig}}
\end{figure}

There are several known bijections between $\mathfrak{S}_n$ and the set of Laguerre 
histories of size $n$. The Foata-Zeilberger bijection $\Psi_{FZ}$ has the property that 
the weight of $\Psi_{FZ}(\sigma)$ is $y^{wex(\sigma)}q^{cr(\sigma)}$. The Françon-Viennot 
bijection $\Psi_{FV}$ has the property that the weight of $\Psi_{FV}(\sigma)$ is 
$y^{asc(\sigma)}q^{31\hbox{-}2(\sigma)}$. All this is present in \cite{Co} and references 
therein, we recall here briefly the definition of $\Psi_{FV}$ since we use it later. Let
$\sigma\in\mathfrak{S}_n$, $j\in\{1,\dots,n\}$ and $k=\sigma(j)$. Then the $k$th step of 
$\Psi_{FV}(\sigma)$ is:
\begin{itemize}
\item a step $\nearrow$ if $k$ is a {\it valley}, {\it i.e.} $\sigma(j-1)>\sigma(j)<\sigma(j+1)$,
\item a step $\searrow$ if $k$ is a {\it peak}, {\it i.e.} $\sigma(j-1)<\sigma(j)>\sigma(j+1)$,
\item a step $\rightarrow$ if $k$ is a {\it double ascent}, {\it i.e.}
$\sigma(j-1)<\sigma(j)<\sigma(j+1)$, or a {\it double descent}, {\it i.e.}
$\sigma(j-1)>\sigma(j)>\sigma(j+1)$.
\end{itemize}
Moreover the weight of the $k$th step is $y^\delta q^i$ where $\delta=1$ if $j$ is an 
ascent and $0$ otherwise, and $i$ is the number of $u\in\{1,\dots,j-2\}$ such that 
$\sigma(u)>\sigma(j)>\sigma(u+1)$ ({\it i.e.} the number of patterns 31-2 such that $j$ 
correspond to the 2). See Figure \ref{fig} for an example.

\medskip

The set of alternating permutations $\mathfrak{A}_n$ consists of permutations 
$\sigma\in\mathfrak{S}_n$ such that $\sigma(2i-1)>\sigma(2i)<\sigma(2i+1)$ 
for any $i\in\{1,\dots,\lfloor \frac n2\rfloor\}$. When $n$ is even, we have
$\sigma\in\mathfrak{A}_n$ if and only if $\sigma$ has no double descent, and no 
double ascent. When $n$ is odd, we have $\sigma\in\mathfrak{A}_n$ if and only if 
$\sigma$ has no double descent, and only one double ascent at position $n$.
Another possible definition is that $\sigma$ is alternating if
$\sigma(1)<\sigma(2)>\sigma(3)<\dots$, but with the convention that $\sigma(0)=0$ and
$\sigma(n+1)=n+1$ the present definition is more natural.

\smallskip

\begin{rem} \label{rem-even}
When $n$ is even, the fact that $\sigma\in\mathfrak{A}_n$ if and only if $\sigma$ has no 
double descent and no double ascent implies the combinatorial interpretation (\ref{combint}).
Indeed, if $\sigma$ is alternating then the corresponding Laguerre history $\Psi_{FV}(\sigma)$
has no horizontal step, and when $y=1$ these Laguerre histories are precisely the weighted Dyck 
paths corresponding to the continued fraction defining $q$-secant numbers. For the case when 
$n$ is odd, see Remark \ref{rem-odd} below.
\end{rem}

\smallskip

\subsection{Eulerian polynomials and their $q$-analogs}
We define
\begin{equation}
A_n(y,q) = \sum_{\sigma\in\mathfrak{S}_n} y^{wex(\sigma)} q^{cr(\sigma)}
\qquad \hbox{and} \qquad
B_n(y,q) = \sum_{\sigma\in\mathfrak{D}_n} y^{wex(\sigma)} q^{cr(\sigma)},
\end{equation}
so that the left-hand sides of (\ref{eq1}) and (\ref{eq2}) are respectively $A_n(-1,q)$ and
$B_n(-\frac1q,q)$.

\begin{prop} \label{incl} We have inversion formulas:
\begin{equation}
   A_n(y,q) = \sum_{k=0}^n \tbinom nk y^{n-k} B_k(y,q) \quad\hbox{and}\quad
   B_n(y,q) = \sum_{k=0}^n \tbinom nk (-y)^{n-k} A_k(y,q).
\end{equation}
\end{prop}

\begin{proof} This is perhaps the most common application of the inclusion-exclusion 
principle. To get this refined version, we have to check how the statistics are changed 
when adding a fixed point to a permutation. The number of crossings does not change, 
and the number of weak exceedances increases by 1 due to the added fixed point. This 
explains the powers of $y$ in the formulas.
\end{proof}

From the formula for $A_n(y,q)$ given in equation (\ref{An}) and the previous proposition,
we obtain that
\begin{equation} \label{Bn}
B_n(y,q)=
 \tfrac 1{(1-q)^n} \sum_{k=0}^n (-1)^k
   \Bigg( \sum_{j=0}^{n-k} y^j  C(n,k,j) \Bigg)
   \left( \sum_{i=0}^k y^iq^{i(k+1-i)}   \right), 
\end{equation}
\[ \hbox{where}\quad C(n,k,j) = \tbinom nj \sum_{i=0}^j q^{j-i}\tbinom ji \tbinom{n-j}{i+k}
    - \tbinom n{j-1} \sum_{i=0}^j q^{j-i}\tbinom {j-1}{i-1} \tbinom{n-j+1}{i+k+1}.
\]
Unfortunately it seems that there is no further simplification, it is just a straightforward
rearrangement of $\sum_{k=0}^n \tbinom nk (-y)^{n-k} A_k(y,q)$. We omit details. Moreover
the two different proofs of (\ref{An}) in \cite{MJV,CJPR} can be modified so as to give a 
direct derivation of (\ref{Bn}).

\smallskip

\subsection{Permutation tableaux} Let $\lambda$ be a Young diagram in English notation,
eventually with empty rows. A {\it permutation tableau} $T$ of shape $\lambda$ is a filling
of $\lambda$ with 0s and 1s such that there is at least a 1 in each column, and for any 0
either all entries to its left are also 0s or all entries above are also 0s. We refer to
\cite{SW} for more details.
We denote by $r(T)$ the number of rows of $T$, by $c(T)$ the number of columns, by
$o(T)$ the number of 1s, by $so(T)=o(T)-c(T)$ the number of {\it superfluous} 1s 
({\it i.e.} all 1s except the topmost of each column).

Via the Steingrímsson-Williams bijection \cite{SW}, the number of fixed points 
(respectively of weak exceedances, of crossings) in a permutation is the number of 
zero-rows (respectively of rows, of superfluous 1s) in the corresponding permutation 
tableau. Let $PT_n$ be the set of permutation tableaux of half-perimeter $n$, and 
$DT_n\subset PT_n$ the subset of permutation tableaux with no zero-row (that we call
{\it derangement tableaux}). Then we have:
\begin{equation} \label{permtab}
A_n(y,q) = \sum_{T\in PT_n} y^{r(T)}q^{so(T)},
\quad\
B_n(y,q) = \sum_{T\in DT_n} y^{r(T)}q^{so(T)},
\end{equation}
\begin{equation} \label{dertab}
 B_n(-\tfrac1q,q) = \sum_{T\in DT_n} (-1)^{r(T)}q^{o(T)-n}.
\end{equation}

\smallskip

\section{Eulerian numbers and $q$-tangent numbers}
\label{sec3}

First we give a characterization of Laguerre histories corresponding to 
permutations $\sigma\in\mathfrak{S}_{n+1}$ such that $\sigma(n+1)=1$.
It is particularly similar to Viennot's criterion \cite[Chapter II, page 24]{XGV}, 
but since the settings are different we believe it is worth writing a direct proof.

\begin{lem} \label{histrest}
Let $\sigma\in\mathfrak{S}_{n+1}$ and $H$ its image by the 
Françon-Viennot bijection as described in Section 2. Then the following conditions 
are equivalent:
\begin{itemize}
\item $\sigma(n+1)=1$
\item Except the first step, $H$ has no step starting at height $h$ with 
  weight $yq^h$.
\end{itemize}
In particular, if this condition is true there is no return to height 0 before the last step.
\end{lem}

\begin{proof}
To begin, suppose that $\sigma(n+1)\neq 1$. Let $k$ be the minimal integer such that 
$\sigma^{-1}(k)>\sigma^{-1}(1)$, and $h$ the starting height of the $k$th step in
$\Psi_{FV}(\sigma)$. We will show that the weight of this $k$th step is $yq^h$.

\medskip

Let $j=\sigma^{-1}(k)$. We have $\sigma(j)<\sigma(j+1)$, otherwise the inequalities
$\sigma(j+1)<\sigma(j)=k$ and $j+1>1$ would imply that the integer $k'=\sigma(j+1)$
contradicts the minimality of $k$. This shows that $j$ is an ascent, so the $k$th step
in $\Psi_{FV}(\sigma)$ is either $\rightarrow$ or $\nearrow$ and has weight $yq^i$
for some $i\in\{0,\dots,h\}$. We have to show that $i=h$. By definition of $\Psi_{FV}$, 
this means that we have to find distinct integers $u_1,\dots,u_h \in \{1,\dots,j-2\}$ such 
that $\sigma(u_l)>\sigma(j)>\sigma(u_l+1)$ for any $l\in\{1,\dots,h\}$.

\medskip

Let $a<k$ be such that the $a$th step of $\Psi_{FV}(\sigma)$ is a step $\nearrow$. 
Let $c=\sigma^{-1}(a)$ and $b$ the minimum integer such that 
$\sigma(b) > \sigma(b+1) > \dots > \sigma(c)$. We have $b<c<j$.
We distinguish two cases.
\begin{itemize}
\item If $\sigma(b)>k$, then there is $u\in\{b,\dots,c-1\}$ such that $\sigma(u)>k>\sigma(u+1)$.
  So we have found one of the $u_l$.
\item Otherwise, we have $\sigma(b-1)<\sigma(b)>\sigma(b+1)$ and $\sigma(b)<k$. So the 
  $\sigma(b)$th step in $\Psi_{FV}(\sigma)$ is a step $\searrow$ before the $k$th step.
\end{itemize}
Thus, to each step $\nearrow$ among the $k-1$ first ones in $\Psi_{FV}(\sigma)$, we can associate
either a $u_l$ satisfying the desired properties, or a step $\searrow$ among the $k-1$ first ones. 
This is sufficient to find indeed the $u_1,\dots,u_h$ and this completes the first part of the proof.
We can use similar arguments to prove the reverse implication.
\end{proof}

\begin{rem} 
It would be possible to consider two other parameters
counting the right-to-left minima and the left-to-right maxima. Indeed, they respectively
correspond to steps starting at height $h$ with weight $yq^h$, and to steps $\rightarrow$ with
weight $yq^0$ or $\searrow$ with weight $q^0$. In this broader context the previous lemma is
immediate because permutations satisfying $\sigma(n+1)=1$ are exactly the ones with only one 
right-to-left minima.
\end{rem}

\smallskip

We now give the proof of Theorem \ref{th1}.

\begin{proof} For any $\sigma\in\mathfrak{S}_n$, we define $\tilde\sigma\in
\mathfrak{S}_{n+1}$ by $\tilde\sigma(i) = \sigma(i)+1$ if $i\in\{1,\dots,n\}$ 
and $\tilde\sigma(n+1)=1$. It is clear that $asc(\sigma)=asc(\tilde\sigma)$ and
$31\hbox{-}2(\sigma)=31\hbox{-}2(\tilde\sigma)$. 

Then for any $\sigma\in\mathfrak{S}_n$, Let $f(\sigma)=\Psi_{FV}(\tilde\sigma)$. 
From the previous lemma, $f$ is a bijection between
$\mathfrak{S}_n$ and the set of Laguerre histories of $n+1$ steps such that:
\begin{itemize}
\item the weight of an horizontal step at height $h$ is $q^i$ or $yq^i$ for some 
  $i\in\left\{0,\dots,h\right\}$
\item Except the first step which has weight $yq^0$, the weight of a step $\nearrow$ 
(resp. $\searrow$) starting at height $h$ is $yq^i$ (resp. $q^i$)
for some $i\in\left\{0,\dots,h-1\right\}$.
\end{itemize}
Moreover the weight of $f(\sigma)$ is $y^{asc(\sigma)}q^{31\hbox{-}2(\sigma)}$. 

In such a Laguerre history, we can remove the first and last step, and
obtain a large Laguerre history with respect to the shifted origin.
In the large Laguerre history, when $y=-1$ the weights on horizontal steps
cancel each other, so that $A_n(-1,q)$ is a sum of weights of Dyck paths of length $n-1$.
When $n$ is even, there is no such Dyck paths so $A_n(-1,q)=0$.

When $n$ is odd, each Dyck path has $\frac {n-1}2$ steps $\nearrow$. So we can factorize 
the sum by $y^{\frac {n-1}2} = (-1)^{\frac {n-1}2}$, and the remaining weighted Dyck 
paths are precisely the ones given by the combinatorial interpretation of the $q$-tangent
numbers. So in this case we have $A_n(-1,q)=(- 1)^{\frac {n-1}2} E_n(q)$.
\end{proof}

\begin{rem} \label{rem-odd}
Through the previous proof we have seen a bijection between permutations and large
Laguerre histories. We easily see that when $n$ is odd, $\sigma\in\mathfrak{S}_n$ is 
alternating if and only if $\tilde\sigma$ is alternating, and this is equivalent to the
fact that $f(\sigma)$ has no horizontal step. When $y=1$ we obtain exactly the
weighted Dyck paths corresponding to the countinued fraction defining the $q$-tangent 
numbers, so this proves the combinatorial interpretation (\ref{combint}) when $n$ is odd.
\end{rem}

\smallskip

Now we prove Theorem \ref{th3}. 

\begin{proof} From Theorem \ref{th1} we have $E_{2n+1}(q)=(-1)^nA_{2n+1}(-1,q)$.
There are interesting cancellations if we directly set $y=-1$ in the proofs of (\ref{An}). 
However, since it is not particularly useful to rewrite these proofs with $y=-1$, we only 
show how to simplify $A_{2n+1}(-1,q)$ and get the right-hand side of (\ref{En}). From 
(\ref{An}) we have:
\begin{equation} \label{An2}
A_{2n+1}(-1,q) = \tfrac {1}{(1-q)^{2n+1}}\sum_{k=0}^{2n+1} (-1)^k\Big( g(n,k) + g(n,k+2) \Big)
\sum_{i=0}^{k} (-1)^iq^{i(k+1-i)},
\end{equation}
where $g(n,k)$ is the sum
\[ g(n,k) = \sum_{j=0}^{2n+1-k}(-1)^j\binom{2n+1}{j}\binom{2n+1}{j+k}. \]
The main step is to simplify this sum. The first term is $\binom{2n+1}k$ and the quotient of 
two successive terms is:
\[
\frac{(-1)^{j+1}\binom{2n+1}{j+1}\binom{2n+1}{j+1+k}}{(-1)^{j}\binom{2n+1}{j}\binom{2n+1}{j+k}}
= - \tfrac{(2n+1-j)(2n+1-k-j)}{(j+1)(j+1+k)} 
= - \tfrac{(-2n-1+j)(-2n-1+k+j)}{(1+j)(k+1+j)}.
\]
so $g(n,k)$ is the hypergeometric series $\binom{2n+1}{k}{}_2F_1( -2n-1+k,-2n-1;k+1;-1)$.
We can use Kummer's summation formula \cite[Chapter 1]{GR}, which reads:
\[ {}_2F_1( a,b;1+a-b;-1)=\frac{\Gamma(1+a-b)\Gamma(1+\frac a2)}
{\Gamma(1+a)\Gamma(-b + 1+\frac a2)}. \]
We cannot rigorously apply the summation formula here because the function
$\Gamma$ is singular at non-positive integers, but we can handle this since
$\Gamma(-m+\epsilon)\sim \frac{(-1)^m}{m!\epsilon}$ for any $m\in\mathbb{N}$
and small $\epsilon$. We write formally:
\[g(n,k) = \frac{\Gamma(2n+2)}{\Gamma(k+1)\Gamma(2n+2-k)}\cdot
\frac{\Gamma(k+1)\Gamma(-n+\frac 12+\frac k2)}
{\Gamma(k-2n)\Gamma(n+\frac 32 +\frac k2)}. \]
When $k$ is even, the only singular value is $\Gamma(k-2n)$ in the denominator 
so $g(n,k)=0$. When $k$ is odd, the singular values compensate:
\[  \lim_{\epsilon \to 0}
 \frac{\Gamma(-n+\frac 12+\frac k2+\epsilon)}{\Gamma(k-2n+2\epsilon)}
 = 2(-1)^{n+\frac{k-1}2}  \frac{\Gamma(2n+1-k)}{\Gamma(n+\frac 12-\frac k2)}
\]
After some simplification it comes that when $n$ is odd,
\[ g(n,k) 
 = \frac{  2(-1)^{n+\frac{k-1}2} \Gamma(2n+2) }
   {(2n+1-k)\Gamma(n+\frac 12-\frac k2)\Gamma(n+\frac 32+\frac k2)}
 = (-1)^{n+\frac{k-1}2} \binom{2n+1}{n-\frac{k-1}2}. \]
Going back to (\ref{An2}), we can restrict the sum to odd $k$ and reindex it so that $k$ 
becomes $2k+1$. Since $A_{2n+1}(-1,q)=(-1)^nE_{2n+1}(q)$ this gives the  formula claimed 
in Theorem \ref{th3}.
\end{proof}

\medskip

\section{Derangements and $q$-secant numbers}
\label{sec4}

First we give the proof of Theorem \ref{th2}.

\begin{proof} The image of a derangement by the Foata-Zeilberger bijection (as defined 
in \cite{Co}) is a Laguerre history such that there is no horizontal step with weight $yq^0$.
This just means that a fixed point cannot be part of a crossing in a permutation. So $B_n(y,q)$ 
is the sum of weights of Motzkin paths of length $n$ such that:
\begin{itemize}
\item the weight of a step $\nearrow$ starting at height $h$ is $y[h+1]$,
\item the weight of a step $\searrow$ starting at height $h$ is $[h]$,
\item the weight of a step $\rightarrow$ starting at height $h$ is $(1+yq)[h]$.
\end{itemize}
When $y=-1/q$, the weights cancel each other on the horizontal steps. So in this case
we can restrict the sum to Dyck paths instead of Motzkin paths, and
the end of the proof is similar to the one of Theorem \ref{th1}.
If $n$ is odd, there is no Dyck path of length $n$ so $B_n(-\frac 1q,q)=0$. 

When $n$ is even, each Dyck path has $\frac n2$ steps $\nearrow$. So we can factorize 
the sum by $y^{\frac n2} = (-\frac 1q)^{\frac n2}$, and the remaining weighted Dyck 
paths are precisely the ones given by the combinatorial interpretation of the $q$-secant 
numbers. So in this case we have $B_n(-\frac 1q,q)=(-\frac 1q)^{\frac n2} E_n(q)$.
\end{proof}

The fact that $B_n(-\frac 1q,q)=0$ when $n$ is odd is also a consequence of 
(\ref{dertab}), because the transposition of derangement tableaux changes the parity
of the number of rows and does not change the number of 1s.

\medskip

To prove Theorem \ref{th4}, it should be possible to use (\ref{Bn}) and specialize
at $y=-\frac 1q$. However it seems more simple to prove it directly with the 
methods and previous results of \cite{MJV}, where the first proof of (\ref{An}) was
given. 

\medskip

A method for computing $A_n(y,q)$ is a Matrix Ansatz. This is a consequence of
the combinatorial interpretation of the PASEP partition function, given by Corteel
and Williams in terms of permutation tableaux \cite{CW}.
There is a similar result to compute $B_n(y,q)$.

\begin{prop} Suppose that we have matrices $D$, $E$, a row vector $\bra{W}$, a column 
vector $\ket{V}$ such that
\begin{equation} \label{ansatz}
DE-qED = I+qE+D, \qquad \bra{W}E=0, \qquad D\ket{V}=0, \qquad \braket{W}{V}=1,
\end{equation}
where $I$ is the identity matrix. Then we have:
\[ B_n(y,q) = \bra{W}(yD+E)^n\ket{V}.\]
\end{prop}

This is to understand as follows. Via the commutation relation, $(yD+E)^n$ can be
uniquely written as a linear combination $\sum_{i,j\geq 0}c_{ij}E^iD^j$, and then 
$B_n(y,q)$ is the constant term $c_{00}$. We give two proofs of this.

\begin{proof} The first proof relies on the combinatorial interpretation of the
3-parameter partition function of the PASEP \cite{CW}. If we have
\[ D'E'-qE'D'=D'+E', \quad \bra{W}E'=0, \quad D'\ket{V}=\ket{V}, \]
then by \cite[Theorem 3.1]{CW}, $\bra{W}(yD'+E')^n\ket{V}$ is the generating function
of permutation tableaux of half-perimeter $n+1$, with no 1 in the first row, where $y$ 
counts the number of rows minus 1, and $q$ counts the number of superfluous 1s. In such 
permutation tableaux, we can remove the first row (which is filled with 0s) and obtain 
any permutation tableaux of half-perimeter $n$. So we get:
\[ A_n(y,q) =  \bra{W}(yD'+E')^n\ket{V}. \]
Together with Proposition \ref{incl}, this gives:
\[ B_n(y,q) =  \bra{W}(yD'-yI+E')^n\ket{V}. \]
But it is straightforward to check that $D=D'-I$ and $E=E'$ satisty conditions (\ref{ansatz}).
\end{proof}
\begin{proof} In the second proof, we sketch how a recursive enumeration of derangement
tableaux directly leads to the matrix Ansatz (\ref{ansatz}). The method is the same as in
\cite[Section 2]{MJV}, or also \cite{AV}. For any word $w$ of size $n$ in $D$ 
and $E$, we define a Young diagram $\lambda(w)$ by the following rule: the South-East 
boundary of $\lambda(w)$ is obtained by reading $w$ from left to right, and drawing a step
East for each letter $D$ in $w$, and drawing a step North for each letter $E$ in $w$.
Let $T_w=\sum_T q^{so(T)}$ where 
the sum is over derangement tableaux of shape $\lambda(w)$. Then we have
\[ T_{w_1DEw_2} =  qT_{w_1EDw_2} + T_{w_1w_2} + qT_{w_1Ew_2} + T_{w_1Dw_2},
\quad\hbox{and}\quad  T_{Ew}=T_{wD}=0.  \]
This is the same kind of relation as the one given by Williams for \Le-diagrams \cite{LW}.
It is obtained by examining a particular corner of the Young diagram, and distinguishing
four different cases whether this corner contains a topmost or non-topmost 1, and a  leftmost
or non-leftmost 1.
We can translate these relations in terms of operators $D$ and $E$ satisfying conditions
(\ref{ansatz}), and such that $T_w=\bra{W}w\ket{V}$. Since $(yD+E)^n$ expands into a (weighted) 
sum of all words in $D$ and $E$ we get all possible shapes for a derangement tableau, 
so that $\bra{W}(yD+E)^n\ket{V}$ is 
the generating function $\sum_wy^{w_D}T_w = B_n(y,q)$, where $w_D$ is the number of letters
$D$ in the word $w$. See the references for more details about the method.
\end{proof}

Now we prove Theorem \ref{th4}.

\begin{proof} We define
\[ \hat D = \tfrac{q-1}{q}\big(D+\tfrac{q}{q-1}\big),
\qquad \hat E = (q-1)\big(E+\tfrac{1}{q-1}\big) \]
so that an immediate computation gives:
\[ \hat D \hat E -q \hat E\hat D = \frac{1-q}q , \qquad \bra{W}\hat E = \bra{W}, 
\qquad \hat D\ket{V}= \ket{V}.\]
and 
\[ (-\tfrac1qD + E)^n = \frac{1}{(q-1)^n}(-\hat D + \hat E)^n.  \]
With the results of \cite{MJV}, it is known that $\bra{W}(-\hat D + \hat E)^n\ket{V}$
is a generating function for weighted involutions and has the following expression:
\[  
\bra{W}(-\hat D + \hat E)^n\ket{V} = \sum_{0\leq k\leq j\leq n} (-1)^j
\left(\tbinom nk - \tbinom n{k-1}\right)q^{(j-k)(n-j-k)-k}.
\]
Since $E_{2n}(q) = (-q)^n B_{2n}(-\tfrac 1q,q)$, we obtain
\[ 
E_{2n}(q)  = \frac 1{(1-q)^{2n}} \sum_{0\leq k\leq j\leq 2n} (-1)^{n+j}
\left(\tbinom {2n}k - \tbinom {2n}{k-1}\right)q^{(j-k)(2n-j-k)+n-k}.
\]
Discarding some negative powers of $q$, we can restrict the sum to $k$ such that $k\leq n$, and then
substitute $k$ with $n-k$, which gives
\[
E_{2n}(q)  = \frac 1{(1-q)^{2n}} \sum_{k=0}^n \sum_{j=n-k}^{2n} (-1)^{n+j}
\left(\tbinom {2n}{n-k} - \tbinom {2n}{n-k-1}\right)q^{(j+k-n)(n+k-j)+k}.
\]
Then we substitute $j$ with $n-k+j$ and obtain
\[
E_{2n}(q)  = \frac 1{(1-q)^{2n}} \sum_{k=0}^n \sum_{j=0}^{n+k} (-1)^{k+j}
\left(\tbinom {2n}{n-k} - \tbinom {2n}{n-k-1}\right)q^{j(2k-j)+k}.
\]
Discarding other negative powers of $q$, we can restrict the sum to $j$ such that $j\leq 2k$ and
this gives the desired expression.
\end{proof}

\begin{rem}
The idea of introducing the operators $\hat D$ and $\hat E$ and their
combinatorial interpretation is at the origin of one of the proof
of (\ref{An}). Similarly, the computation above can be refined to
obtain another proof of (\ref{Bn}).
\end{rem}

\smallskip

\section{Second proof of the $q$-secant and $q$-tangent formulas}
\label{sec5}

A second proof of (\ref{An}), inspired by Penaud's bijective proof of the Touchard-Riordan
formula \cite{JGP}, was given in \cite{CJPR}. Actually, a slight modification of the
method can give a direct proof of (\ref{Bn}). In this section, we show that this can
be used to obtain Theorem \ref{th3} and \ref{th4} directly from the definition of $E_n(q)$
in terms of weighted Dyck paths. We begin with the $q$-secant numbers because this case is 
more simple.

\subsection{The $q$-secant numbers}
First, we notice that $(1-q)^{2n}E_{2n}(q)$ is the sum of weights of Dyck paths $H$
of length $2n$ such that:
\begin{itemize}
\item the weight of a step $\nearrow$ starting at height $h$ is 1 or $-q^{h+1}$,
\item the weight of a step $\searrow$ starting at height $h$ is 1 or $-q^h$.
\end{itemize}

\begin{prop} \label{decomp}
There is a weight-preserving bijection between these Dyck paths $H$ and
couples $(H_1,H_2)$ such that for some $k\in\{0,\dots,n\}$,
\begin{itemize}
\item $H_1$ is a left factor of a Dyck path, of $2n$ steps and final height $2k$,
\item $H_2$ is a weighted Dyck path of length $2k$, with the same weights as $H$
 and also the condition that there is no two consecutive steps $\nearrow\searrow$ 
 both weighted by 1.
\end{itemize}
\end{prop}
\begin{proof} This is a direct adaptation of \cite[Lemma 1]{CJPR} so we only sketch
the proof. The idea is to look for the maximal factors of $H$ being Dyck paths and 
having only steps with weight 1. To obtain $H_1$ from $H$ we transform into a step
$\nearrow$ any step which is not inside one of these maximal factors. And to obtain 
$H_2$ from $H$ we just remove these maximal factors.
\end{proof}

By an elementary recurrence, the number of left factors of Dyck paths of $2n$ steps and final 
height $2k$ is $\binom {2n}{n-k} - \binom {2n}{n-k-1}$.
So to obtain Theorem \ref{th4} it remains only to prove the following proposition.

\begin{prop} \label{core1}
The sum of weights of Dyck paths of length $2k$, satisfying conditions
as $H_2$ above, is
\[  M_k(q) = \sum_{j=0}^{2k}(-1)^{j+k}q^{j(2k-j)+k}. \]
\end{prop}

\begin{proof} It is possible to adapt the proof of \cite[Proposition 5]{CJPR}.
However we take a slightly different point of view here, and show that we can 
exploit some properties of T-fractions. The combinatorial theory for T-fractions
was developed by Roblet and Viennot in \cite{RV}, and gives generating functions
for the kind of paths we are dealing with here. Indeed, this reference gives 
immediately that:
\begin{equation}  \label{fracM}
  \sum_{k=0}^{\infty} M_k(q) t^k =  
  \cfrac{1}{1+t-\cfrac{(1-q)^2t}{1+t-\cfrac{(1-q^2)^2t}{\ddots}}}.
\end{equation}
To make this article more self-contained, we briefly show how to obtain (\ref{fracM}).
Let us consider weighted Schröder paths, such that:
\begin{itemize}
\item the weight of a step $\nearrow$ starting at height $h$ is 1 or $-q^{h+1}$,
\item the weight of a step $\searrow$ starting at height $h$ is 1 or $-q^{h}$,
\item the weight of a step $\longrightarrow$ is $-1$.
\end{itemize}
The sum of weights of these Schröder paths of length $2k$ is $M_k(q)$. Indeed there is
a sign-reversing involution on Schröder paths, such that the fixed points are the Dyck 
paths with no two consecutive steps $\nearrow\searrow$ both weighted by 1. The involution
is obtained by exchanging the first occurrence of $\longrightarrow$ with weight $-1$, or
$\nearrow\searrow$ with weight 1, with (respectively) $\nearrow\searrow$ with weight 1, or
$\longrightarrow$ with weight $-1$. For the generating function of Schröder paths, standard
arguments \cite{Fla} gives the continued fraction in (\ref{fracM}).

\smallskip

We need to extract the coefficients of $t^k$ in the continued fraction. To do this, there 
are well-established methods linking this kind of continued fractions with quotients of
contiguous basic hypergeometric series \cite[Chapter 19]{CU}. In the present case, we can 
use a limit case of (19.2.11a) in \cite{CU}. If we consider the more general continued 
fraction:
\[ M(z) = \cfrac{1}{1+t-\cfrac{(1-qz)^2t}{1+t-\cfrac{(1-q^2z)^2t}{\ddots}}}, \]
then \cite[(19.2.11a)]{CU} with $c=tq$ and $a=-b=i\sqrt{qt}$ gives that 
\[ M(z) = \frac{1}{1-z} \cdot 
   \frac{ {}_2\phi_1(a,b;ab;q;qz)  } {  {}_2\phi_1(a,b;ab;q,z) }.
\]
It is the same series as in 
\cite{CJPR} but with different values for $a$ and $b$, so we proceed similarly
with a Heine's transformation \cite[Chapter 1]{GR}:
\[ {}_2\phi_1(a,b;ab;q,z) = \frac{(az,b;q)_\infty}{(ab,z;q)_\infty} 
   {}_2\phi_1(a,z;az;q,b).\]
So
\[ M(z) = 
  \frac{(aqz,b,ab,z;q)_\infty}{ (1-z) (ab,qz,az,b;q)_\infty} \cdot
  \frac{ {}_2\phi_1(a,qz;aqz;q,b)  } {  {}_2\phi_1(a,z;az;q,b) }
= \frac{1}{1-az} \cdot
  \frac{ {}_2\phi_1(a,qz;aqz;q,b)  } {  {}_2\phi_1(a,z;az;q,b) },
\]
and
\begin{eqnarray*}
 M(1)  & = & \frac{1}{1-a} {}_2\phi_1(a,q;aq;q,b) 
   = \sum_{n=0}^\infty \frac{b^n}{1-aq^n} \\
    &  =  & \sum_{n=0}^\infty\sum_{j=0}^\infty b^na^jq^{jn} = 
\sum_{n=0}^\infty\sum_{j=0}^\infty t^\frac{n+j}{2}  (-i)^n i^j  q^{jn+\frac{n+j}2}. 
\end{eqnarray*}
To obtain the coefficient of $t^k$ in $M(1)$, we just have to restrict the sum to 
the couples $(n,j)$ such that $n=2k-j$. This gives the desired formula.
\end{proof}

\subsection{The $q$-tangent numbers}
To prove Theorem \ref{th4} with the same method, we first need to modify slightly 
our expression for $E_{2n+1}(q)$. Let $P_k$ be the polynomial 
$\sum_{j=0}^{2k+1}(-1)^{j+k}q^{j(2k+2-j)}$. By elementary properties of the binomial 
coefficients we have
\[
E_{2n+1}(q) = \frac 
 {\sum\limits_{k=0}^n\left(\tbinom {2n+1}{n-k} - \tbinom {2n+1}{n-k-1}\right) P_k} 
 {(1-q)^{2n+1}}
 = \frac{ \sum\limits_{k=0}^n\left(\tbinom {2n}{n-k} - \tbinom {2n}{n-k-1}\right)
 \tfrac {P_k+P_{k-1}}{1-q} } {(1-q)^{2n}}.
\]
The latter expression is the one we can prove with the previous method.
We decompose the weighted paths exactly as in the case of $q$-secant numbers,
{\it i.e.} as in Proposition \ref{decomp}. The sole difference is treated in the
following proposition, which is therefore the last step of the second proof of 
Theorem \ref{th4}.

\begin{prop} \label{core2}
The sum of weights of Dyck paths of length $2k$, with weights $1$ or
$-q^{h+1}$ for any step ($\nearrow$ or $\searrow$) starting at height $h$,
and no two consecutive steps $\nearrow\searrow$ both weighted by 1, is
\[  N_k(q) = \frac {P_k+P_{k-1}}{1-q}.\]
\end{prop}

\begin{proof} We follow the scheme of Proposition \ref{core1}. The Schröder paths
in this case give the following T-fraction:
\begin{equation}  \label{fracN}
  \sum_{k=0}^{\infty} N_k(q) t^k =  
  \cfrac{1}{1+t-\cfrac{(1-q)(1-q^2)t}{1+t-\cfrac{(1-q^2)(1-q^3)t}{\ddots}}}.
\end{equation}
We need to consider the more general continued fraction:
\[ N(z) = \cfrac{1}{1+t-\cfrac{(1-qz)(1-q^2z)t}{1+t-\cfrac{(1-q^2z)(1-q^3z)t}{\ddots}}}.\]
And then \cite[(19.2.11a)]{CU} with $c=tq$ and $a=-b=i\sqrt{t}q$ gives that 
\[ N(z) = \frac{1}{1-z} \cdot 
   \frac{ {}_2\phi_1(a,b;tq;q;qz)  } {  {}_2\phi_1(a,b;tq;q,z) }.
\]
In this case the Heine's transformation gives:
\[  {}_2\phi_1(a,b;tq;q,z) =
   \frac{ (-i\sqrt{t}q, i\sqrt{t}qz ;q)_\infty  }{ (z,tq;q)_\infty }
 {}_2\phi_1(i\sqrt{t}, z ; i\sqrt{t}qz; q; -i\sqrt{t}q)
\]
Hence, after some simplification:
\[
N(1) =  \frac{1}{1-i\sqrt{t}q} {}_2\phi_1(i\sqrt{t}, q ; i\sqrt{t}q^2; q; -i\sqrt{t}q  )
  =
\sum_{n=0}^\infty \frac{ (1-i\sqrt{t}) (-i\sqrt{t}q)^n }
{( 1-i\sqrt{t}q^n )( 1-i\sqrt{t}q^{n+1} )}.
\]
The usual method to reduce degrees in denominators leads to:
\[
(1-q)N(1)=-\sum_{n=0}^\infty \frac {1-i\sqrt{t}}{1-i\sqrt{t}q^n}(-i\sqrt{t})^{n-1}
          +\sum_{n=0}^\infty \frac {1-i\sqrt{t}}{1-i\sqrt{t}q^{n+1}}(-i\sqrt{t})^{n-1}.
\]
Some terms in this sum cancel each other and it remains:
\[
(1-q)N(1)=- \sum_{n=0}^\infty \frac 1{1-i\sqrt{t}q^n}(-i\sqrt{t})^n
          + \sum_{n=0}^\infty \frac 1{1-i\sqrt{t}q^{n+1}}(-i\sqrt{t})^{n-1} 
          + \frac{(i\sqrt{t})^{-1}}{1-i\sqrt{t}}.
\]
At this point we can expand the fractions as in the previous case, and
we readily obtain that the coefficient of $t^k$ in $(1-q)N(1)$ is $P_k+P_{k-1}$.
\end{proof}

\section{Final remarks}

It might be possible to give a combinatorial proof of Proposition \ref{core1} or \ref{core2},
with a sign-reversing involution whose fixed points account for the terms in the right-hand 
side (as was done for the Touchard-Riordan formula in \cite{JGP}). In either case, we
could not even conjecture a set of possible fixed points for such an involution.

\bigskip

To our knowledge there is no simple parity-independent expression of the classical 
Euler numbers $E_n=E_n(1)$. So it is remarkable that Theorems \ref{th3} and \ref{th4}
are so similar. One might ask if there is a parity-independent formula for $E_n(q)$.
Actually from Theorems \ref{th3} and \ref{th4} we get that for any even or odd $n$
\[
 (-1)^{\frac{n-1}2}A_n(-1,q) = \tfrac{(-1)^{\lfloor n/2 \rfloor}}{(1-q)^n} 
 \sum_{k=0}^{\lfloor \frac n2 \rfloor}
 \left( \tbinom{n}{k} - \tbinom{n}{k-1}  \right)
 \sum_{i=0}^{n-2k} (-1)^{k+i} q^{i(n-2k-i)+i}, 
\] 
(the right-hand side is 0 when $n$ is even because 
$\sum_{i=0}^{n-2k}(-1)^{i+k}q^{i(n+1-2k-i)}=(-1)^{n+1-k}$ and
$\sum_{k=0}^{n}(-1)^k(\tbinom {2n}k - \tbinom {2n}{k-1}) =0 $
), and also 
\[
 (-1)^{\frac n2}B_n(-\tfrac 1q,q) = \tfrac{(-1)^{\lfloor n/2 \rfloor}}{(1-q)^n} 
 \sum_{k=0}^{\lfloor \frac n2 \rfloor}
 \left( \tbinom{n}{k} - \tbinom{n}{k-1}  \right)
 \sum_{i=0}^{n-2k} (-1)^{k+i} q^{i(n-2k-i) + \frac n2 -k} , 
\] 
(the right-hand side is 0 when $n$ is odd because $\sum_{i=0}^{n-2k}(-1)^{i+k}q^{i(n-2k-i)}=0$ ),
and adding the previous two identities gives:
\[
 E_n(q) = \tfrac{(-1)^{\lfloor n/2 \rfloor}}{(1-q)^n} 
 \sum_{k=0}^{\lfloor \frac n2 \rfloor}
 \left( \tbinom{n}{k} - \tbinom{n}{k-1}  \right)
 \sum_{i=0}^{n-2k} (-1)^{k+i} q^{i(n-2k-i)}( q^i+q^{\frac n2 -k} ).
\]


\bigskip

\end{document}